\documentclass[11pt]{amsart}
\usepackage{amssymb}
\usepackage[]{latexsym}
\usepackage{amsthm}

\newtheorem{thm}{Theorem}
\newtheorem{cor}[thm]{Corollary}
\newtheorem{lem}[thm]{Lemma}
\newtheorem{slem}[thm]{Sublemma}

\newtheorem{prop}[thm]{Proposition}

\newtheorem*{mr}{Main result}

\theoremstyle{remark}
\newtheorem{rem}[thm]{Remark}
\newtheorem{exmp}[thm]{Example}

\theoremstyle{definition}

\newenvironment{pf}{\par\noindent{\bf Proof.}\enspace\ignorespaces}{\qed\par\par}
\def\qed{\relax\ifmmode\hskip2em \Box\else\unskip\nobreak\hskip1em $\Box$\fi}

\newcommand{\cha}{\mbox{char}}

\newcommand{\bQ}{{\mathbb{Q}}}
\newcommand{\bF}{{\mathbb{F}}}

\newcommand{\bP}{{\mathbb{P}}}
\newcommand{\bZ}{{\mathbb{Z}}}
\newcommand{\bC}{{\mathbb{C}}}

\begin{document}

\title{Squares in arithmetic progression over number fields}
\author{Xavier Xarles} \email{xarles@mat.uab.cat}
\address{Departament de Matem\`atiques\\Universitat
Aut\`onoma de Barcelona\\08193 Bellaterra, Barcelona, Catalonia}

\thanks{partially supported by the grant MTM2006--11391.}

\begin{abstract}
We show that there exists an upper bound for the number of squares
in arithmetic progression over a number field that depends only on
the degree of the field. We show that this bound is 5 for
quadratic fields, and also that the result generalizes to
$k$-powers for $k>1$.
\end{abstract}

\maketitle

In this note we are dealing with the following natural problem:
Given a number field $K$, is there a maximum for the number of
distinct elements $a_0, \dots, a_n$ in $K$ such that
$a_i^2-a_{i-1}^2=a_{i+1}^2-a_i^2$ for $i=1, \dots , n-1$? We will
prove that there is a bound for this maximum, and that this bound
only depends on the degree of the field $K$ over $\bQ$. In fact we
will show that the same result is also valid for $k$-powers, now
the bound depending also on $k$.

The problem has a long history for $K=\bQ$. In a letter written to
Frenicle in 1640, Fermat proposed the problem of proving that there
are no four squares in arithmetic progression. Euler gave the first
published proof of this result in 1780. In a different direction, in
1970 Szemer{\'e}di proved that there exists at most $o(N)$ squares in an
arithmetic progression of length $N$, and this result was improved
by Bombieri, Granville and Pintz in 1992 \cite{BGP} to
$O(N^{2/3}(\log N)^A)$ for a suitable constant $A$ studying the
arithmetic progressions that contain 5 squares, and by Bombieri and
Zannier in 2002 in \cite{BZ} to $O(N^{3/5}(\log N)^A)$ for a
suitable constant $A$ studying the ones that contain 4 squares.

The question for higher powers has also a long history. It is
known that it does not exists a nontrivial three term arithmetic
progression of $k$-th powers for $k \ge 3$. Observe that, when $k$
is odd, we do have non constant three term arithmetic progression
of $k$-th powers, the ones of the form $-a^k$, $0$ and $a^k$ for
$a\in \bQ$. In these cases, for non-trivial three term arithmetic
progression we mean non constant and with $a_1\ne 0$. The cases
$k=3$ and $k=4$ are mentioned in Carmichael's 1908 book on
diophantine equations. The cases $k=5, \dots 31$ were done by
Denes in 1952 \cite{De}. The cases that $k\ge 17$ is a prime
number congruent to 1 modulo 4 where done by Ribet \cite{Ri}, and
the rest of the cases by Darmon and Merel in 1997 \cite{DM}.

The problem is related to some concrete curves having only trivial
rational points (trivial in some sense). The rational points of
this curves determine arithmetic progressions having squares at
the first $n$ terms, and the trivial points correspond to the
constant arithmetic progression.

For example, four consecutive squares in an arithmetic progression
give a rational point in an elliptic curve, that one can show has
only 8 solutions, all coming from the constant arithmetic
progression.

The main result of this note is the following theorem.

\begin{mr}
For any $d\ge 1$, there exists a constant $S(d)$ depending only
$d$ such that, if $K/\bQ$ verifies that $[K:\bQ]=d$ and
$a_i:=a+i\;r$ is an arithmetic progression with $a$ and $r \in K$,
and $a_i$ are squares in $K$ for $i=0,1,2,\cdots,S(d)$, then $r=0$
(i.e. $a_i$ is constant).

Furthemore, if $d=2$, then $S(2)=6$.
\end{mr}

In some forthcoming papers it is studied over which quadratic
fields we have fourth squares in arithmetic progression ( by E.
Gonz{\'a}lez-Jim{\'e}nez and J. Steuding \cite{GJS}), and five squares in
arithmetic progression (by E. Gonz{\'a}lez-Jim{\'e}nez and X. Xarles
\cite{GX}).

This note is organized as follows. In the first section we
translated the problem to some problem concerning the
determination of all the rational points of some algebraic curves
$C_n$, and we prove some preliminary results. The second section
we give a lower bound for the gonality of these curves, which we
use in section 3 in order to obtain the existence of the constant
$S(d)$. In section 4 we investigate the value $S(2)$, proving some
results concerning the rational points of $C_4$ and $C_5$ over
quadratic fields. Finally, in the last section we show how to
proof the result on $k$-powers, and we comment some
generalizations of the problem.

This paper had its origin in a question asked by Ignacio Larrosa
Ca{\~n}estro, which I could not answer satisfactorily (see the last
section). I thank him, and also Adolfo Quir{\'o}s, Enric Nart, Joaquim
Ro{\'e}, Javier Cilleruelo and Henry Darmon  for some conversations
and comments on the subject. I especially thank Andrew Granville
and Qing Liu for some fruitful comments and for giving me some
useful references and Enrique Gonz{\'a}lez Jim{\'e}nez for some comments
and corrections.

\section{Translation to algebraic curves}

We say that some elements $a_0,\dots,a_n$ on a field $K$ are in
arithmetic progression if there exists $a$ and $r$ elements of
$K$, $a\cdot r\ne 0$, such that $a_i=a+i\;r$ for any
$i=0,\dots,n$. This is equivalent, of course, of having
$a_i-a_{i-1}=r$ fixed for any $i=1,\dots,n$.

First of all, observe that, in order to study squares inb
arithmetic progressions, we can and will identify the arithmetic
progressions $\{a_i\}$ and $\{a'_i\}$ such that there exists a
$c\in K^*$ with $a'_i=c^2a_i$ for any $i$. In case that $K=\bQ$,
we can suppose then that $a$ and $r$ are coprime integers.

Now, suppose that there exists three squares $x_0^2$, $x_1^2$ and
$x_2^2$ in arithmetic progression over $K$. This means that
$$x_1^2-x_0^2=x_2^2-x_1^2$$
and so $(x_0,x_1,x_2)$ is a solution of the equation
$$f(X_0,X_1,X_2):=X_0^2-2X_1^2+X_2^2=0.$$ We are not interested in the trivial
solution $(0,0,0)$, and solutions that are equal up to
multiplication by an element in $K^*$ (in fact, in $(K^*)^2$) we
consider them equal. So we will work with solutions of the
projective curve $f(X_0,X_1,X_2)$ in the projective plane $\bP^2$.

Similarly, in order to consider $n+1$ squares in arithmetic
progression, with $n\ge 2$, we will take the curve $C_{n}$ in
$\bP^n$ determined by the $n-1$ equations
$$ f(X_i,X_{i+1},X_{i+2})=0 \ \mbox{ for } i=0,\dots, n-2.$$
The assignment of the arithmetic progression corresponding to any
point induces the following map $\varphi_{n}\colon C_{n} \to
\bP^1$ given by $\varphi_n(X_0,\dots,X_n)=[X_0^2:X_1^2-X_0^2]$.

The curve $C_{n}$ is a non-singular projective curve over any
field $K$ of characteristic bigger that $n$, as we will prove in
the following lemma.

\begin{lem}\label{reduction} Let $n\ge 1$ and let $K$ be any field,
$p=\cha(K)$, with $p>n$ or $p=0$. Then the curve $C_{n}$ is a
non-singular projective curve of genus $g_{n}:=(n-3)2^{n-2}+1$.

Moreover, let $\varphi_n \colon C_{n} \to \bP^1$ given by
$\varphi_{n}(X_0,\dots,X_n)=[X_0^2:X_1^2-X_0^2]$. Then $\varphi_n$
has degree $2^n$, and it is ramified at the points above $[i:1]$ for
$i=0,\dots, n$.
\end{lem}

\begin{pf} We use first the jacobian criterium in order to show
non singularity. The Jacobian matrix of the system of equations
defining $C_{n}$ is
$$A:=( \partial f(X_i,X_{i+1},X_{i+2})/\partial X_j)_{0\le i\le n-2\,,\,0\le j\le n}.$$

For any $j_1<j_2$, denote by $A_{j_1,j_2}$ the matrix obtained by
$A$ by deleting the columns $j_1$ and $j_2$; it is an square
matrix of size $(n-1)\times (n-1)$. It is easily shown that its
determinant verifies that
$$|A_{j_1,j_2}|= \pm 2^{n-1} \left(\prod_{i\ne j_1,j_2} X_i \right)
(j_2-j_1).$$

We want to show that, for any $[x_0:\dots:x_n]\in C(K)$, there
exists $\{j_1,j_2\}$ such that $|A_{j_1,j_2}|(x_0:\dots:x_n)\ne
0$.

The first crucial observation is that any point
$[x_0:\dots:x_n]\in C(K)$ can have at most one $i=1,\dots,n$ such
that $x_i=0$. If the characteristic of the field is zero this is
clear. If the characteristic is $p$, suppose that $x_i=x_j=0$ with
$i<j$. Since $x_i^2=a+ir$ for certain $a$ and $r\in K$, we will
have that $(i-j)r=0$, so $r=0$, which implies $a=0$, which is not
possible, or $i-j=0$ in $K$, so $i+kp=j$ for certain $k\in \bZ$,
which again is not possible if $p>n$.

So, if $p>n$ or $p=0$, for any point $[x_0:\dots:x_n]\in C(K)$, if
all $x_i$ are different from 0 then $|A_{j_1,j_2}|\ne 0$ $\forall
j_1\ne j_2$, and if $x_i=0$, then  $|A_{i,j_1}|\ne 0$ $\forall
j_1\ne i$, hence the rank of $A$ is $n-1$.

The genus of the curve $C_{n}$ can be computed by induction on $n$
applying the Hurwitz formulae to the natural forgetful cover of
degree $2$
$$C_{n} \to C_{n-1},$$
which is ramified on the $2$ points with $x_n=0$. Or can be
computed by using the map $\varphi_n\colon C_{n} \to \bP^1$ which
has degree $2^n$. The ramification points are the points
$[x_0:\dots:x_n]$ such that there exists some $i$ with $x_i=0$.
Such points have image by $\varphi_n$ equal to $[i:1]$, and have
ramification index equal to $2$.
\end{pf}

We will call the $2^n$ points $[\pm 1:\dots:\pm 1]$ the trivial
points. They correspond to the points $P$ such that
$\varphi_{n}(P)=[0:1]=\infty$, so giving the constant arithmetic
progression. So the first aim of this note is to prove that for $n$
sufficiently large with respect to $d$, they are the only
$K$-rational points for any extension $K/\bQ$ of degree $d$.

Firstly, one can easily prove the existence of such a bound but
depending on the field $K$. Observe that, if $n>3$, then the genus
is bigger than $1$, so, by Faltings' Theorem (previously known as
the Mordell Conjecture), for any number field $K$ and $n>3$,
$C_n(K)$ is finite. One can prove even more.

\begin{lem}\label{Kbound} Let $K/\bQ$ be a finite extension. Then there exists a
constant $n_K$ such that $C_{n_K}(K)$ has only the trivial points.
\end{lem}
\begin{pf} Consider first the finite set of points in $\varphi_4(C_4(K))\subset
\bP^1(K)$. We want to show that $\varphi_n(C_n(K))\subset
\varphi_4(C_4(K))$ is equal to $\{\infty\}$ for $n$ sufficiently
large. This is equivalent to show that for any $P:=[a:r]\in
\varphi_4(C(K))$ not equal to $\infty$, there exists some $n_P$
such that $P$ is not in $\varphi_n(C_n(K))$. But this is obvious
from the following sublemma.\end{pf}

\begin{slem} Let $K/\bQ$ be a finite extension, and let $\{a_i\}$
a non-constant arithmetic progression. Then there exists some $n$
such that $a_n$ is not a square in $K$.
\end{slem}
\begin{pf} We can reduce to the case that all $a_n$ are in the
ring of integers in $K$. Consider a prime ideal $\mathfrak p$ in
the ring of integers of $K$ with residue field a prime field
$\bF_p$, $p>2$, and such that the sequence $\{\widetilde{a_i}\}$,
reduction modulo $\mathfrak p$ of $\{a_i\}$, is not constant. Then
the sequence $\{\widetilde{a_i}\}$ can have only $(p+1)/2$
consecutive squares, and hence also $\{a_i\}$. So there exists
some $n\le (p+3)/2$ such that $a_n$ is not a square in $K$.
\end{pf}

\begin{rem} The constant $n_K$ in the lemma is not constructive since
it depends on being able to find all the $K$-rational points of
$C_n$ for some $n>3$.
\end{rem}

So we have proved that there exists a bound depending only on the
field $K$. In order to show that one can find a bound depending
only on the degree of the field, we will apply a criterium of Frey
that is a consequence of Faltings' Theorem. To do this we need to
give a lower bound for the gonality of the curves $C_n$.

\section{The gonality of $C_n$ over $\bQ$.}

Recall that the gonality $\gamma(C_K)$ of a curve $C$ over a field
$K$ is the minimum $m$ such that there exists a morphism
$\phi:C\to \bP^1$ of degree $m$ defined over $K$. For example,
hyperelliptic curves have gonality $2$. The aim of this section is
to give a lower bound for the gonality $\gamma_n$ of $C_n$ over
$\bQ$.

From the forgetful map $C_n\to C_{n-1}$, and using that $C_0$ has
genus 0, and hence it has gonality 1, we get the upper bound
$2^{n-2}\ge \gamma_n$. In order to get a lower bound, we will use
the following result which is well known for the experts (see, for
example, Proposition 3 in \cite{Fre}).

\begin{prop}\label{rgonality} Let $C$ be a curve over a number
field, and let $\wp$ be a prime of good reduction of the curve,
with residue field $\bF_q$. Denote also by $C'$ the reduction of
the curve $C$ modulo $\wp$. Then the gonality $\gamma$ of $C$
verifies that
$$\gamma \ge \frac{\sharp C'(\bF_{q^n})}{q^n+1}$$ for any $n\ge 1$.
\end{prop}

\begin{pf}  Suppose there is
a map $f:C\to \bP^1$ defined over $K$ of degree $\gamma$. First,
we want to show that there is a morphism $f':C'\to \bP^1$ of
degree $\gamma'\le \gamma$. This result is well-known, and can be
deduced from results by Abhyankar (see for example \cite{Na}), and
also from Deuring \cite{Deu} or from Lemma 5.1 in \cite{N-S}, but
we write a short proof for completion.

I learned these proof form Q. Liu. By applying Lemma 4.14 in
\cite{LL}, with $X=C$, $Y=\bP^1$, and $\mathcal X$ a smooth model
of X, we get a model $\mathcal Y$ together with a rational map
$\mathcal X \to \mathcal Y$ which is quasi-finite on the (good)
reduction $C$ of $X$. Considering the image $D$ of $C'$ in the
reduction of $\mathcal Y$, then $D$ is rational curve (because the
reduction of $\mathcal Y$ has arithmetical genus 0). The degree of
$C' \to D$ is less than the total degree $\gamma$. Now the
normalization $D'$ of $D$ is smooth rational and geometrically
irreducible because $C'$ is geometrically irreducible. Since the
curve $C$ is a conic over a finite field, we get a rational point
on $D'$ (because of Waring's result), hence $D' = \bP^1$.

Consider this morphism $f':C\to \bP^1$ defined over $\bF_q$ of
degree $\gamma'\le \gamma$. Hence, at most $\gamma'$ elements in
$C(\bF_{q^n})$ can go to the same point in $\bP^1(\bF_{q^n})$,
which have $q^n+1$ points. So $(q^n+1)\gamma \ge(q^n+1)\gamma'\ge
\sharp C(\bF_{q^n})$.
\end{pf}

\begin{cor}\label{gonality} For any $n\ge 3$, let $\gamma_n$ be the gonality of
$C_n$. Then $2^{n-2}\ge\gamma_n\ge 2^{n-1}/n$.
\end{cor}

\begin{pf} Let $p$ be a prime such that $2n>p>n$. Then $C_n$ has
good reduction over $p$, and its reduction is given by the same
curve $C_n$. Consider the $2^n$ trivial points
$\varphi_n^{-1}(\infty) \subset C_n(\bF_p)$. By applying the
proposition we get that
$$\gamma \ge \frac{\sharp C_{n}(\bF_{p})}{p+1} \ge \frac {2^n}{p+1}\ge \frac{2^{n}}{2n}$$
by using that $p\le 2n-1$.
\end{pf}

\begin{rem} The result in the corollary can be improved for some concrete
values of $n$ by considering the some prime such that $p>n$ and
some $m$ (usually $m=1$ or $2$), and consider also all the points
in $C_n(\bF_{p^m})$. On the other hand, one can also use some
results concerning the gonality over $\bC$, and using that
$\gamma(C_{\bQ})\ge \gamma(C_{\bC})$.
\end{rem}

\section{Proof of the main theorem, first part}

First of all, observe that in order to show the existence of a
constant $S(d)$ such that for any degree $d$ extension $K/\bQ$, the
only arithmetic progressions with $S(n)$ consecutive squares is the
constant ones, we only need to show that the existence of such a
constant such that $C_n(K)$ contains only the trivial points for
$n=S(d)$.

We will use the following criterium of Frey \cite{Fre}, proved
also by Abramovich in his thesis.

\begin{thm}[(Frey)] \label{Frey}
Let $C$ a curve over a number field $K$, with gonality $\gamma>1$
over $K$. Fix an algebraic closure $\overline{K}$ of $K$ and
consider the points of degree $d$ of $C$,
$$C^d(K):=\bigcup_{[L:K]\le d} C(L) \ \subset C(\overline{K})$$
where the union is over all the finite extensions of $K$ inside
$\overline{K}$ of degree $\le d$. Suppose that $2d<\gamma$. Then
$C^d(K)$ is finite.
\end{thm}

Hence, by proposition \ref{gonality}, we get that, if $2d<
2^{n-1}/n<\gamma_n$, then there exists a finite number of points
in $C_n^d(\bQ)$. So, there exists only a finite number of
extensions $K_i/\bQ$ of degree $d$ such that $C(K_i)$ contains
some non trivial point, and for any other $K/\bQ$ of degree $d$,
$C(K)$ contains only the trivial points.

For any such extension, we apply now the lemma \ref{Kbound}, so
there exists some constant $n_{K_i}$ such that $C_{n_{K_i}}$
contains only the trivial points. Hence, considering
$S(d):=\max_{i}\{n_{K_i}\}$ we get the result.

\begin{rem} The constant $S(d)$ as defined above is not explicit
since it depends on Faltings' theorem twice, one for the
construction of the extensions $K_i$, and, second, in the
computation of the constants $n_{K_i}$, as remarked before.

However, one could guess that the constant $S(d)$ will not be much
bigger than the one verifying $d< 2^{S(d)-2}/S(d)$, so the correct
value of $S(d)$ should be $O(\log(d))$.
\end{rem}

\section{Proof of main theorem, second part: the case $d=2$}

In this section we want to calculate explicitly the constant $S(2)$.
Observe that we do have 5 squares in arithmetic progression over
quadratic fields, for example $$7^2=49, 13^2=169, 17^2=289, 409
,23^2=529$$ over $\bQ(\sqrt{409})$, so $S(2)>5$. In fact, we have an
infinite number of examples of such progressions, even we do have
only a finite number over any fixed field.

\begin{lem}\label{4infinity} There is an infinite number of different quadratic extensions
$K/\bQ$ such that $C_4(K)$ contains non trivial points.\end{lem}

\begin{pf} Consider the curve parametrizing the arithmetic
progressions $a_i$ such that $a_i$ is an square for $i=0,1,2,4$.
This curve is a genus 1 curve given by equations
$$X_0^2-2X_1^2+X_2^2=0 \ \mbox{ and }\ 3X_0^2-4X_1^2+X_4^2=0.$$
Since it has a (trivial) point, it is isomorphic to an elliptic
curve, that can be given by the Weierstrass equation
$$y^2=x(x+2)(x+6)$$
One shows by standard methods that this curve has an infinite
number of rational points, and in fact $E(\bQ)\cong \bZ \oplus
\bZ/2\bZ$, being $(2,8)$ a generator of the torsion free part.

Now, for any point $P\in E(\bQ)$, we consider the associated
progression $a_0:=x_0^2$, $a_1:=x_1^2$, and so on. By considering
the field $K:=\bQ(\sqrt{a_3})$, we get that in $K$ the arithmetic
progression $a_i$ has 5 consecutive squares. Equivalently,
consider the degree 2 map $\psi:C_4\to E$, and the points $Q\in
C_4(K)$ such that $\psi(Q)=P$.  Since by Faltings theorem there
are only a finite number of points in $C_4(K)$, and we have a
infinitely many points in $E(\bQ)$, we get infinitely many such
fields $K$.
\end{pf}

In a forthcoming paper \cite{GX}, we study over which quadratic
fields one has 5 squares in arithmetic progression, an we get, for
example, that $\bQ(\sqrt{409})$ is the smallest (in terms of the
discriminant) of such fields.

Observe that the gonality $\gamma$ of $C_5$ over $\bQ$ is bounded
below by $2^4/5$ by using corollary \ref{gonality}, but in fact
one can show the gonality is strictly bigger that 4. So by Frey
result, Theorem \ref{Frey}, we get that there is only a finite
number of quadratic fields having non trivial points. We will show
that in fact there is none.

In order to show this, we will study in more detail the case of 5
squares. We will prove that, if $K/\bQ$ is a quadratic extension and
$P\in C_4(K)$, then $\phi(P)\in \bP^1(\bQ)$, hence the progression
is defined over $\bQ$ (as they are all the ones obtained from the
lemma \ref{4infinity}). This result will imply that $C_5(K)$ contain
only the trivial points.

To show the result on $C_4$ we need first to study how are the
points of $C_3(\bQ)$ and $C_3(K)$ for $K/\bQ$ of degree 2. Observe
that $C_3$ is isomorphic to an elliptic curve once fixed a rational
point, and we will take always $[1:1:1:1]$. We get then a group
operation $\oplus$ on $C_3$. The following lemma describes some easy
cases.

\begin{lem} For any field $K$, consider $P:=[x_0: x_1: x_2: x_3]\in C_3(K)$,
 and let $Q=[\pm 1:\pm 1: \pm 1: \pm 1]$ be any trivial point. Then
$\ominus P=[x_3:x_2:x_1:x_0]$ and
$$\ominus P\oplus Q=\left\{ \begin{array}{rl}
\mbox{Case 1: } &  [\pm x_0: \pm x_1: \pm x_2: \pm x_3], \\
\mbox{Case 2: } &  [\pm x_3: \pm x_2: \pm x_1: \pm x_0],
\end{array}
\right.$$ where the Case 1 is when the number of $-1$ in $Q$ is
even, and the signs of the $x_i$'s are the same that the signs of
the corresponding coordinate of $Q$; and the Case 2 is when the
number of $-1$ in $Q$ is odd, with the same rule for the signs.
Moreover, in the case 1 the point $Q$ has order 2, and in Case 2 has
order 4.
\end{lem}

\begin{pf} This is elementary, and we will show only some cases. Moreover,
it can be shown using two distinct strategies.

The first one is to work out the formulae for the addition on an
elliptic curve given as intersection of two quadrics in $\bP^3$,
together with a rational point $O$. We explain the general
procedure. Consider the plane that passes through $O$ with
multiplicity three or four; there is only one, which is called the
osculating plane. This plane intersects the curve in another point
$O'$, the osculation point, that is the same $O$ if the
multiplicity of intersection is equal to four. In order to compute
the sum of two points $P_1$ and $P_2$, consider the plane that
passes through $P_1$, $P_2$ and $O'$ (in case $P_1=P_2$, it should
cut the curve in $P_1$ with multiplicity two). Denote the fourth
intersection point by $P_3'$. Finally, consider the plane that
passes through $O$, $O'$ and $P_3'$. The other intersection point
is the point $P_3$. Observe that what this procedure states is
that three (different) points $P_1$, $P_2$ and $P_3$ add to 0 if
and only if all three together with $O'$ are coplanar. If we want
to compute $-P$ for a given point $P$, we consider the plane that
passes through $P$, $O$ and $O'$; the fourth intersection point is
$-P$, as is easily seen from the procedure for the sum.

In our case, $O=[1:1:1:1]$, the osculating plane is given by the
equation $-X_0+3X_1-3X_2+X_3$ and the osculating point by
$O':=[1:-1:-1:1]$. Then the symmetric of a point $P=[x_0: x_1 :
x_2: x_3]$ is $\ominus P=[x_3:x_2:x_1:x_0]$, since $P$, $\ominus
P$, $O$ and $O'$ are coplanar.

Consider now the point $Q=[-1:1:1:1]$. For any given point
$P=[x_0: x_1: x_2: x_3]$, let $P':=[x_0: x_1: x_2: -x_3]$. Then
$P$, $P'$, $0'$ and $Q$ are coplanar, so $P\oplus Q\oplus P'=0$.
This means that $P\oplus Q=-P'=[-x_3:x_2:x_1:x_0]$. Using the same
argument one shows that
$$P\oplus [-1:1:1:1]=[-x_3:x_2:x_1:x_0],$$
$$P\oplus [1:-1:1:1]=[x_3:-x_2:x_1:x_0], $$
$$P\oplus [1:1:-1:1]=[x_3:x_2:-x_1:x_0],$$
and
$$P\oplus [1:1:1:-1]=[x_3:-x_2:x_1:-x_0]. $$

The other cases are easily obtained from these. For example, one
has that $[-1:1:1:1]\oplus [1:-1:1:1]=[-1:1:-1:1]$, so
$$P\oplus [-1:1:-1:1]=P\oplus[-1:1:1:1]\oplus
[1:-1:1:1]=$$$$=[-x_3:x_2:x_1:x_0]\oplus
[1:-1:1:1]=[x_0:-x_1:x_2:-x_3].$$

Observe that the operation $\oplus [-1:1:1:1]$ has order 4, hence
the point $[-1:1:1:1]$ has order 4. One gets that a trivial point
has order 2 if it has an even number of $-1$, and order 4
otherwise. Observe that $Q+Q=[-1:1:1:-1]$ in all order 4 cases.

The second approach uses the following observation: since $C_3$
has no CM (a fact that can be shown just by having a non integral
$j$-invariant, which is $2^413^3/3^2$), the only automorphisms of
order 2 of $C_3$ (as a genus 1 curve) which fix the 0 are the
identity automorphism and negation, and the other automorphisms of
order 2 are the translations by a 2 torsion point and negation
followed by translations. We get then that
$$[x_0: x_1: x_2: x_3] \mapsto [x_3:x_2:x_1:x_0]$$
is the negation, since it fixes $O$ and its not the identity. Once
we know this, we get that $\ominus Q=Q$ if and only if $Q$ is a
trivial point with an even number of $-1$, hence in this case $Q$
has order 2. But
$$[x_0: x_1: x_2: x_3] \mapsto [\pm x_0: \pm x_1: \pm x_2: \pm
x_3]$$ has order 2, and imposing that $O+Q=Q$ and $Q+Q=Q$, we get
the correct choice of signs, solving the Case 1. Now, for the case
2, consider first the point $Q=[-1:1:1:1]$. One has that $\ominus
Q=\ominus [-1:1:1:1]=[1:1:1:-1]$, and that $[-1:1:1:1]\oplus
[1:-1:-1:1]=[1:1:1:-1]$, so $2*[-1:1:1:1]=[1:-1:-1:1]$ and hence $Q$
has order 4. From this fact one easily deduce that all the other 2
points have order 4. The formulae are then easily obtained once we
know it for some specific $P$.
\end{pf}

Now we need to recall the classical result of Fermat, concerning
squares in arithmetic progressions over $\bQ$, together with some
other cases.

\begin{lem}\label{Fermat} Let $F$ be one of the following genus 1 curves, given
as intersection of two quadrics in $\bP^3$:
$$C_3: X_0^2-2X_1^2+X_2^2=0 \ \mbox{ and }\ X_1^2-2X_2^2+X_3^2=0,$$
$$F_1: 2X_0^2-3X_1^2+X_2^2=0\  \mbox{ and }\ X_1^2-3X_2^2+2X_3^2=0, $$
$$F_2: X_0^2-3X_1^2+2X_2^2=0\  \mbox{ and }\ 2X_1^2-3X_2^2+X_3^2=0, $$
Then $F(\bQ)=\{\pm 1:\pm 1: \pm 1: \pm 1]\}$.

As a consequence, if $\{a_i\}$ is an arithmetic progression over
$\bQ$ such that $a_0$, $a_1$, $a_2$ and $a_3$ are squares, or
$a_0$, $a_1$, $a_3$ and $a_4$ are squares, or $a_0$, $a_2$, $a_3$
and $a_5$ are squares, then $\{a_i\}$ is constant.
\end{lem}

\begin{pf} We consider all the cases
at the same time. Being $F$ a genus one curve with some rational
point, for example $[1,1,1,1]$, $F$ is isomorphic to its jacobian.
One shows by standard methods that their jacobian are the elliptic
curves given respectively by the Weierstrass equations
$y^2=x(x-1)(x+3)$, $y^2=x(x-1)(x+8)$ and $y^2=x(x-4)(x+5)$.
These elliptic curves have only a 8 rational
points, as proved by standard descent methods. Hence $F(\bQ)$ has
only 8 points, which must be $\{[\pm1:\pm1:\pm1:\pm1]\}$.

The assertion about the arithmetic progressions is easy now. The
first case we already considered, so we consider the second one.
Suppose $a_0=x_0^2$, $a_1=x_1^2$, $a_3=x_2^2$ and $a_4=x_3^2$ are
squares. Then $2(x_1^2-x_0^2)=x_2^2-x_1^2=2(x_3^2-x_2^2)$, getting
the two equations defining $F_1$. Since the only rational
solutions of $F_1$ correspond to having $a_i=1$ for all $i$, the
progression is constant. The case of $F_2$ is similar.
\end{pf}

\begin{cor} Let $K/\bQ$ a degree 2 extension, and let $\sigma$ the
generator of the Galois group. Let $P=[x_0: x_1: x_2: x_3] \in
C_3(K)$ be a non trivial point. Then $$P^{\sigma}:=[\sigma(x_0):
\sigma(x_1): \sigma(x_2): \sigma(x_3)]=\left\{ \begin{array}{rl}
\mbox{Case 1:} &  [\pm x_3: \pm x_2: \pm x_1: \pm x_0]\\
\mbox{Case 2:} &  [\pm x_0: \pm x_1: \pm x_2: \pm x_3].
\end{array}
\right.$$ Furthermore, in the case 2, $\phi_3(P)\in \bP^1(\bQ)$.
\end{cor}

\begin{pf} Let $P=[x_0: x_1: x_2: x_3] \in
C_3(K)$, and consider $Q:=P\oplus P^{\sigma}$. Since
$Q^{\sigma}=Q$, $Q\in C_3(\bQ)$, so $Q=[\pm 1:\pm 1: \pm 1: \pm
1]$. Hence $P\ominus Q=\ominus P^{\sigma}$, and by the lemma above
we have that either $P\ominus Q= [\pm x_0: \pm x_1: \pm x_2: \pm
x_3]$ in the case 1, hence $P^{\sigma}=[\pm x_3: \pm x_2: \pm x_1:
\pm x_0]$, or $P\ominus Q= [\pm x_3: \pm x_2: \pm x_1: \pm x_0]$;
hence $P^{\sigma}=[\pm x_0: \pm x_1: \pm x_2: \pm x_3]$. In this
last case we have that $\sigma(x_i^2)=x_i^2$ (may be after
rescaling the coordinates), and, therefore, the corresponding
arithmetic progression is defined over $\bQ$.
\end{pf}

\begin{exmp} We is easy to give examples of the case 2, for example
the one given by
$x_0=1$, $x_1=5$, $x_2=7$ and $x_3=\sqrt{61}$, in
$K=\bQ(\sqrt{61})$. There are also examples of the case 1: take
$K=\bQ(\sqrt{13})$, and consider $x_0=1$, $x_1=10+3\sqrt{13}$,
$x_2=15+4\sqrt{13}$ and $x_3=18+5\sqrt{13}$. Then
$P^{\sigma}=[-x_3:x_2:-x_1:x_0]$.
\end{exmp}

Now we apply this results to study the curve $C_4$. There are two
different forgetful maps from $C_4$ to $C_3$, forgetting the first
term and forgetting the last term. We will use this assertion in
order to show the following result.

\begin{prop} Let $K/\bQ$ a degree 2 extension, and let
$P \in C_4(K)$ be a non trivial point. Then $\phi_4(P)\in
\bP^1(\bQ)$. \end{prop}

\begin{pf} Let $P=[x_0: x_1: x_2: x_3: x_4]\in C_4(K)$, that
we can suppose with $x_i\ne 0$ $\forall i=0,\dots,4$, and consider
$P_0:=[x_0: x_1: x_2: x_3]$ and $P_1=[ x_1: x_2: x_3:x_4]\in
C_3(K)$. Suppose that $P_i$ for some $i=0,1$ is in case 2 of the
previous corollary. Then, may be after rescaling the coordinates,
we have that $x_1^2$, $x_2^2$ and $x_3^2$ are in $\bQ$, hence the
arithmetic progression is defined over $\bQ$, which is equivalent
to $\phi_4(P)$ being in $\bP^1(\bQ)$.

So we can suppose both $P_i$ are in case 1. We can take then that
$x_2=1$. Then there exists $\lambda_0$ and $\lambda_1$ in
$K\setminus \{0\}$ such that $\sigma(x_0)=\pm \lambda_0 x_3$,
$\sigma(x_1)=\pm \lambda_0 =\pm \lambda_1 x_4$, $1=\sigma(1)=\pm
\lambda_0 x_1=\pm \lambda_1 x_3$, $\sigma(x_3)=\pm \lambda_0
x_0=\pm \lambda_1$ and $\sigma(x_4)=\pm \lambda_1 x_1$. Hence we
get that $x_1=\pm 1/\lambda_0$, $x_3=\pm 1/\lambda_1$, $
\sigma(x_0)=\pm \lambda_0/\lambda_1$ and $\sigma(x_4)=\pm
\lambda_1/\lambda_0$.  Using now that the $x_i$ belong to $C_4$ we
get easily that $\lambda_0^2+\lambda_1^2=2$ and
$1/\lambda_0^2+1/\lambda_1^2=2$. From these one easily deduce that
$\lambda_i=\pm 1$ or $\pm 2$ for $i=0$ and $1$. But this implies
that all the $x_i=\pm 1$, hence $P$ is already defined over $\bQ$.
\end{pf}

The content of the last proposition is that all the rational
points of $C_4$ defined over quadratic extensions are obtained
taking square roots of some arithmetic progressions over $\bQ$,
and essentially with the method explained in the proof of lema
\ref{4infinity}.

Now, we are going to use this result to show the non existence of 6
squares in arithmetic progressions over quadratic fields.

\begin{thm} Let $K/\bQ$ a degree 2 extension, and suppose
$P \in C_5(K)$. Then $P=[\pm1:\pm1:\pm1:\pm1:\pm1:\pm1]$.
\end{thm}

\begin{pf} let $D$ a square free integer, and consider
$K:=\bQ(\sqrt{D})$. Suppose we have a point $P=[x_0: x_1: x_2: x_3:
x_4:x_5]\in C_5(K)$. By the previous proposition, we have that
$x_i^2\in \bQ$ for all $i=0,\dots,5$. So we can and will suppose
that $x_i^2=a_i:=a+i\;r$ for some $a$ and $r\in \bZ$ coprime
integers. If all the $x_i$ are in $\bQ$, then we are done by
Fermat's result. So suppose we have some $x_i\not{\in}\bQ$, and
hence $x_i^2=Dy_i^2$ for some $y_i\in \bQ$.

Observe first that no integer $D>5$ can divide two of the $x_i^
2$'s; since if a prime $p>5$ divides $x_{i}^2-x_{j}^2=(i-j)r$, then
it divides $(i-j)<6$, which is not possible, or it divides $r$, and
hence it divides $a=x_i^2-i\;r$, which again is not posible being
$a$ and $r$ coprime.

So, if $D>5$, we must have that $x_i\bQ$ for all $i=0,\dots,5$,
except may be for some $i=j$. If $j=0$, $1$, $4$ or $5$, then we
will have 4 squares in arithmetic progression over $\bQ$, so
$x_i^2=1$ for all $i$ by Fermat (or lemma \ref{Fermat}), giving
$x_j^2=1$, contrary to the hypothesis. Now, if $j=2$, then we will
have $a_0$, $a_1$, $a_3$ and $a_4$ are squares over $\bQ$, so again
$x_i^2=1$ by lemma \ref{Fermat}. The same argument shows the case
$j=3$.

Hence we are reduced to the cases $D=2,$ $3,$ $4$ and $5$. The case
$D=4$ is obviously trivial, being $4$ an square. The case $D=5$ is
again easy, since if $5$ divides two elements of the progression,
they must be $a_0$ and $a_5$, and then $a_1,\dots,a_4$ will be 4
squares in arithmetic progression over $\bQ$. For the case $D=3$,
$D$ dividing two elements of the progression, we have that, or
$a_0$, $a_1$, $a_3$ and $a_4$ are squares over $\bQ$, a case that we
already considered, or $a_0$, $a_2$, $a_3$ and $a_5$  are squares
over $\bQ$, which is the third case considered in lemma
\ref{Fermat}, or, finally, $a_1$, $a_2$, $a_4$ and $a_5$ are
squares, which is equivalent to the first case.

So we are reduced to the case $D=2$. First, suppose that
$a_0=2y_0^2$, $a_2=2y_2^2$ and $a_3=y_3^2$ for some $y_i\in \bQ$.
Then we have that $y_0^2+y_3^2=3y_2^2$, which has no solutions over
$\bQ$. This implies that we cannot have that $2$ divides $a_i$ and
$a_{i+2}$ for any $i=0,\cdots,3$. Second, suppose we have
$a_0=2y_0^2$, $a_1=y_1^2$ and $a_3=y_3^2$. Then we have
$4y_0^2+y_3^2=3y_1^2$, which again has no solutions over $\bQ$. With
these fact we solve the remaining cases, proving that $2$ cannot
divide two terms of the $a_i=x_i^2$.
\end{pf}

\section{Some generalizations and conjectures}

In order to compute explicitly the constant $S(d)$ for $d>2$ one
cannot use the same argument we did in the last section. In fact,
we do not even know if $S(3)>4$, since we don't know a way to
produce 5-terms arithmetic progressions of squares over cubic
fields. One possible idea could be to use the natural maps from
$C_n$ to sufficiently many elliptic curves with rank 0, and then
showing that is a formal immersion at $p$ for some prime $p>2$, as
it is done (in a different context) in \cite{Fre}. But this result
will not be sufficient to conclude, because this curves contain
always non trivial points over finite fields with cardinal a power
of $p^{2}$, since all the elements of $\bF_p$ are squares over
such fields. We do not know any other general argument to show
this type of results, which essentially are the computation of all
the rational points in the symmetric product of some curve $C$.

As we mention in the introduction, one can ask also for higher
powers the same question we did for squares. And, in fact, the
same type of arguments work for solving the existence of a uniform
bound. So we get the following result.

\begin{thm}\label{kpow} Given $d\ge 1$ and $k\ge 2$, there exists a constant
$S(d,k)$ depending only on $d$ and $k$ such that, if $K/\bQ$
verifies that $[K:\bQ]=d$ and $a_i:=a+i\;r$ is an arithmetic
progression with $a$ and $r \in K$, and $a_i$ are $k$-powers in
$K$ for $i=0,1,2,\cdots,S(d,k)$, then $r=0$ (i.e. $a_i$ is
constant).
\end{thm}

\begin{pf} First, consider the case $d=1$. In this case, as we
mention in the introduction, it is known after the work of Denes
\cite{De}, Ribet \cite{Ri}, and Darmon and Merel \cite{DM}, that
the only three term arithmetic progression of $k$-th powers for $k
\ge 3$ over $\bQ$ are the constant ones if $k$ is even, and the
constant plus the ones of the form $-a^k$, $0$ and $a^k$ for $a\in
\bQ$ if $k$ is odd. Using these it is clear that there are no non
constant fourth term arithmetic progression of $k$-th powers for
$k$ odd, and hence $S(1,k)=3$ if $k$ even and $=4$ if $k$ odd.

Now, fix $k>1$ and $d>1$. The assertion of the theorem is
equivalent to showing that for $n$ large enough in terms of $d$,
the only points of degree $d$ of the curves $C_{n,k}$ in
$\bP^n_{\bQ}$ determined by the $n-1$ equations
$$ f_k(X_i,X_{i+1},X_{i+2})=0 \ \mbox{ for } i=0,\dots, n-2,$$
where $f_k(X,Y,Z):=X^k-2Y^k+Z^k$ are the trivial points $[\pm
1:\pm 1:\dots:\pm 1]$ if $k$ is even, or $[1:1:\dots: 1]$ is $k$
is odd. The same arguments used to show lemma \ref{reduction}
shows that $C_{n,k}$ has good reduction at any prime $p>n$ and not
dividing $k$. Using the same argument than in section 3, we only
need to show that the gonality of the curves $C_{n,k}$ tends to
infinity when $n$ goes to infinity, so Theorem \ref{Frey} applies
again.

If $k$ is even, the exact same argument than in corollary
\ref{gonality} concerning the lower bound of the gonality applies
if $n>k$ (in order to avoid the primes $p$ dividing $k$), so we
get the result.

If $k$ is odd, the curve $C_{n,k}$ contains only one trivial point
in the reduction, unless the field $\bF_p$ contains some $k$-roots
of unity. In order to have these, we consider primes $p$ of good
reduction such that $p \equiv 1 \pmod{k}$. In this case we have
that $C_{n,k}(\bF_p)$ contains $k^n$ points. Now, we apply a
well-known consequence of Dirichlet's theorem on primes in
arithmetic progressions asserting that there exists a constant
$c(k)$ depending on $k$ such that, for any $n>c(k)$, there exists
a prime $p$ verifying that $n<p<2n$ and $p\equiv 1 \pmod{k}$.

Now the argument works as follows: choose an $n>c(k)$ and $>k$ and
a prime $n<p<2n$ with $p\equiv 1 \pmod{k}$. Then the gonality
$\gamma_{n,k}$ of $C_{n,k}$ verifies, by proposition
\ref{rgonality}, that
 $$\gamma_{n,k} \ge \frac{\sharp C_{n,k}(\bF_p)}{p+1} \ge \frac
 {k^n}{2n}.$$
 Hence, for $n$ large enough with respect to $d$ we get
 $\gamma_{n,k}>2d$, so $C_{n,k}$ contains only a finite number of
points of degree $d$ over $\bQ$, and the proof is finish by using
an analog of lemma \ref{Kbound}.
\end{pf}

In fact, one can ask even more general questions, not just
restricting to powers, but, for example, to images of a fixed
polynomials. In a paper in preparation \cite{XT}, we study in a
even more general context when similar questions have sense and
when can be asked affirmatively, and also in which cases one can
give explicit uniform upper bounds.

On the other hand, once looking for the computation of specific
bounds, like the ones for $d=1$ or $d=2$, one cannot use the same
type of reasoning we used. As we already mention, the case $d=1$
is known, getting $S(1,k)=2$, but the proof uses for most cases
the same type of ideas that for proving Fermat's last theorem,
i.e. Wiles and others ideas. However, we think that the case $k=3$
and $d=2$ could be solved by a similar argument we did for the
case $k=2$ and $d=2$.

A. Granville observed in a personal communication that one can use
the main theorem to prove that there are always $o_d(N)$ squares
in any arithmetic progression over any number field of degree $d$,
just as a simple application of Szem{\'e}redi's theorem on arithmetic
progressions, as he did for the case $d=1$ using Fermat's result.
More generally, we have the following result.

\begin{cor} Let $d\ge 1$ and $k\ge 2$ be integers. Then
there are $o_{d,k}(N)$ $k$-powers of a field of degree $d$ in any
arithmetic progression $a+im$ with $m\ne0$ and $1\le i\le
N$.\end{cor}

\begin{pf} Let $I$ be the set of $i$ for which $a+im$ is a $k$-power.
Then $I$ does not contain any $S(d,k)$-term arithmetic
progressions by the theorem \ref{kpow}. Hence $|I|=o_{d,k}(N)$ by
Szem{\'e}redi's theorem, which states for any $\delta>0$, if $N$ is
sufficiently large, then any subset of $\{1,2,..,N\}$ of size $M>
\delta N$ has an $M$-term arithmetic progression.\end{pf}

This last result can probably be improve to some bound of the type
$O_{d,k}(N^{1-c_{d,k}})$, for some constant $0<c_{d,k}\le
\frac{k-1}{k}$, by using similar arguments of the ones used in
\cite{BGP} and \cite{BZ}. One could even guess if $c_{d,k}$ can
always be taken equal to $\frac{k-1}{k}$, which will generalize
Rudin's conjecture which asserts this in the case $k=2$ and $d=1$.

Finally, let us mention that one can ask a 2-dimensional (or even
$s$-dimensional) analogous question. By these we mean the
following natural question (which it was ask to me, in a different
form, by Ignacio Larrosa Ca{\~n}estro and it was the origin of this
paper): is there a constant $S$ such that, the only degree 2
polynomials $f(x)$ with rational coefficients verifying that
$f(i)$ is a square for $i=0,\dots,S$ are the squares of a (degree
1) polynomials? This question can be translated to the computation
of all the rational points of some algebraic surfaces. Such
question is investigated for example in \cite{Al}, an it seems
that the natural guess is $S=8$. But a positive answer of the
existence of such $S$ is related to the so called Lang's (and
Bombieri) conjecture, which asserts that for a general type
surface the rational points are all contained in a finite number
of strict subvarieties, in our case curves and points. In fact,
for $S=8$ one gets that the corresponding surface is of general
type. But the only known cases of this conjecture, corresponding
to subvarieties of abelian varieties, does not apply in these
cases, being all these surfaces complete intersections in $\bP^n$.

\end{document}